\documentclass[reqno,11pt]{amsart}
\usepackage{amscd,amssymb}

\theoremstyle{plain}
\newtheorem{Thm}[subsection]{Theorem}
\newtheorem{Cor}[subsection]{Corollary}
\newtheorem{Lem}[subsection]{Lemma}
\newtheorem{Prop}[subsection]{Proposition}
\newtheorem{Conj}[subsection]{Conjecture}

\theoremstyle{definition}
\newtheorem{Def}[subsection]{Definition}

\theoremstyle{remark}

\newtheorem{Rem}[subsection]{Remark}

\newtheorem{conj}{Conjecture}

\numberwithin{equation}{section}
\renewcommand{\rm}{\normalshape}

\newif\ifShowLabels
\ShowLabelstrue
\newdimen\theight
\def\TeXref#1{%
    \leavevmode\vadjust{\setbox0=\hbox{{\tt
        \quad\quad  {\small \rm #1}}}%
    \theight=\ht0
    \advance\theight by \lineskip
    \kern -\theight \vbox to
    \theight{\rightline{\rlap{\box0}}%
    \vss}%
    }}%

\ShowLabelsfalse% comment this out if labels should be printed

\renewcommand{\sec}[2]{\section{#2}\label{S:#1}%
    \ifShowLabels \TeXref{{S:#1}} \fi}
\newcommand{\ssec}[2]{\subsection{#2}\label{SS:#1}%
    \ifShowLabels \TeXref{{SS:#1}} \fi}

\newcommand{\refss}[1]{Section ~\ref{SS:#1}}

\newcommand{\reft}[1]{Theorem ~\ref{T:#1}}
\newcommand{\refl}[1]{Lemma ~\ref{L:#1}}
\newcommand{\refp}[1]{Proposition ~\ref{P:#1}}
\newcommand{\refc}[1]{Corollary ~\ref{C:#1}}

\newcommand{\refe}[1]{\eqref{E:#1}}

\newenvironment{thm}[1]%
    { \begin{Thm} \label{T:#1}  \ifShowLabels \TeXref{T:#1} \fi }%
    { \end{Thm} }

\renewcommand{\th}[1]{\begin{thm}{#1} \sl }
\renewcommand{\eth}{\end{thm} }

\newenvironment{lemma}[1]%
    { \begin{Lem} \label{L:#1}  \ifShowLabels \TeXref{L:#1} \fi }%
    { \end{Lem} }
\newcommand{\lem}[1]{\begin{lemma}{#1} \sl}
\newcommand{\elem}{\end{lemma}}

\newenvironment{propos}[1]%
    { \begin{Prop} \label{P:#1}  \ifShowLabels \TeXref{P:#1} \fi }%
    { \end{Prop} }
\newcommand{\prop}[1]{\begin{propos}{#1}\sl }
\newcommand{\eprop}{\end{propos}}

\newenvironment{corol}[1]%
    { \begin{Cor} \label{C:#1}  \ifShowLabels \TeXref{C:#1} \fi }%
    { \end{Cor} }
\newcommand{\cor}[1]{\begin{corol}{#1} \sl }
\newcommand{\ecor}{\end{corol}}

\newenvironment{defeni}[1]%
    { \begin{Def} \label{D:#1}  \ifShowLabels \TeXref{D:#1} \fi }%
    { \end{Def} }
\newcommand{\defe}[1]{\begin{defeni}{#1} \sl }
\newcommand{\edefe}{\end{defeni}}

\newenvironment{remark}[1]%
    { \begin{Rem} \label{R:#1}  \ifShowLabels \TeXref{R:#1} \fi }%
    { \end{Rem} }
\newcommand{\rem}[1]{\begin{remark}{#1}}
\newcommand{\erem}{\end{remark}}

\newenvironment{conjec}[1]%
    { \begin{Conj} \label{Co:#1}  \ifShowLabels \TeXref{Co:#1} \fi }%
    { \end{Conj} }
\renewcommand{\conj}[1]{\begin{conjec}{#1} \sl }
\newcommand{\econj}{\end{conjec}}

\newcommand{\eq}[1]%
    { \ifShowLabels \TeXref{E:#1} \fi
       \begin{equation} \label{E:#1} }
\newcommand{\eeq}{ \end{equation} }

\newcommand{\prf}{ \begin{proof} }
\newcommand{\epr}{ \end{proof} }

\newcommand\del{\delta}

\newcommand\lam{\lambda}		\newcommand\Lam{\Lambda}
\newcommand\sig{\sigma}		

		\newcommand\Ome{\Omega}

\newcommand\calC{{\mathcal{C}}}
\newcommand\calD{{\mathcal{D}}}
\newcommand\calE{{\mathcal{E}}}

\newcommand\calH{{\mathcal{H}}}

\newcommand\calK{{\mathcal{K}}}

\newcommand\calM{{\mathcal{M}}}

\newcommand\calO{{\mathcal{O}}}
\newcommand\calP{{\mathcal{P}}}

\newcommand\calZ{{\mathcal{Z}}}

		\newcommand\bfG{{\mathbf G}}

		\newcommand\bfM{{\mathbf M}}

		\newcommand\bfP{{\mathbf P}}

\newcommand\QQ{\mathbb{Q}}

\newcommand\ZZ{\mathbb{Z}}

\newcommand\CC{\mathbb{C}}

\newcommand\sdp{\times \hskip -0.3em {\raise 0.3ex
\hbox{$\scriptscriptstyle |$}}} % semidirect product

\newcommand\End{\operatorname{End\,}}
\newcommand\Ext{\operatorname{Ext}}

\newcommand\Hom{\operatorname {Hom}}

\newcommand\Id{\operatorname{Id}}

\newcommand\Ind{\operatorname{Ind}}

\newcommand\supp{\operatorname{supp}}

\newcommand\Tr{\operatorname{Tr}}

\newcommand\oB{{\overline{B}}}

\newcommand\oh{{\overline{h}}}

\newcommand\oP{{\overline{P}}}

\newcommand\oPhi{{\overline{\Phi}}}

\newcommand\tilh{{\widetilde{h}}}

\newcommand\ten{\otimes}

\newcommand{\ra}{\rangle}
\newcommand{\la}{\langle}

\newcommand\Irr{\operatorname{Irr}}

\newcommand\ch{\operatorname{ch}}
\newcommand\ocalH{\overline{\calH}}
\newcommand{\oone}{\bar{1}}
\renewcommand{\oPhi}{\overline{\Phi}}
\newcommand\oi{\bar{i}}
\newcommand\ovr{\overline{r}}

\begin{document}
Dedicated to J.~Bernstein on the occasion of his 70th birthday
\title{Bernstein components via Bernstein center}
\author{Alexander Braverman and David Kazhdan (with an appendix by R.~Bezrukavnikov)}
\begin{abstract}
Let $G$ be a reductive $p$-adic group. Let $\Phi$ be an invariant distribution on $G$ lying in the Bernstein center $\calZ(G)$.
We prove that $\Phi$ is supported on compact elements in $G$ if and only if it defines a constant function on every
component of the set $\Irr(G)$; in particular, we show that the space of all elements of $\calZ(G)$ supported on compact elements is
a subalgebra of $\calZ(G)$. Our proof is a slight modification of the argument from Section 2 of \cite{dat}, where our result is proven in one direction.
\end{abstract}
\maketitle

\sec{}{Introduction}

\ssec{}{Components of $\Irr(G)$}
In this paper $G$ denotes the set of point of a connected reductive algebraic group over a local non-archimedian
field $\calK$. We shall denote by $\calM(G)$ the category of smooth complex representations of $G$. This category is equivalent to the category of unital modules over the Hecke algebra $\calH(G)$. We let $\Irr(G)$ denote the set of isomorphism classes irreducible objects of $\calM(G)$.
Bernstein and Zelevinsky defined a decomposition
of the set $\Irr(G)$ of irreducible representations of $G$ into a union of certain components
$\Ome$; this decomposition in fact defines a decomposition of $\calM(G)$ into a product of the corresponding categories. The set $\calC(G)$ of components of $\Irr(G)$ is in one-to-one correspondence with pairs
$(M,\sigma)$ where $M$ is a Levi subgroup of $G$ and $\sigma$ is a cuspidal representation of $M$ (the data of $(M,\sig)$ is uniquely
determined by $\Ome$ up to natural equivalence relation
generated by conjugation and multiplying $\sig$ by an unramified character of $M$). An element $\pi\in \Irr(G)$ lies in $\Ome(M,\sig)$ if and only if there exists a parabolic subrgoup $P$ containing $M$ as a Levi subgroup and an unramified character $\chi$ of $M$ such that $\pi$ is a subquotient
of the induced representation $\Ind_P^G(\sig\ten \chi)$ (here we use the natural map $P\to M$ in order to view $\sig\ten \chi$ as a representation of $P$). Every $\Ome(M,\sig)$ is equipped with a map to an irreducible affine algebraic variety $\overline{\Ome}(M,\sig)$.  The variety $\overline{\Ome}(M,\sig)$ is in fact a quotient of the torus of unramified characters of $M$ by a finite group. The above map has finite fibers and is generically one-to-one. We shall say that a function $f:\Ome(M,\sig)\to \CC$ is regular if it comes from a regular function on $\overline{\Ome}(M,\sig)$. We shall say that a function $f:\Irr(G)\to \CC$ is regular iff it is regular when restricted to every component.

\ssec{}{Bernstein center}
Let $\calZ(G)$ be the center of the category $\calM(G)$. It is easy to see that it consists of all invariant distributions $\Phi$ on $G$ such that
for any $h\in \calH(G)$ we have $\Phi\star h\in \calH(G)$. It is enough to test the above condition for all $h=e_K$ where $K$ is an open compact subgroup of $K$ and $e_K$ is the Haar measure on it.

By Schur-Quillen lemma any $\Phi\in \calZ(G)$ defines a function on the set $\Irr(G)$.
Bernstein proved (cf. \cite{BD}) that in this way we get an isomorphism between $\calZ(G)$ and the algebra of regular functions on $\Irr(G)$.
Thus $\calZ(G)$ has both "geometric" (in terms of distributions on $G$) and "spectral" (in terms of functions on $\Irr(G)$) description. The relationship between these two descriptions tends to be quite non-trivial.
This note is devoted to one particular aspect of this relationship. Namely, we are going to prove the following

\th{main}Let $\calZ_{comp}(G)$ denote the subspace of $\calZ(G)$ consisting of distributions supported on compact elements. Similarly, let
$\calZ_{lc}(G)$ denote the subalgebra of $\calZ(G)$ consisting of those elements $\Phi$ for which $f(\Phi)$ is a locally constant function (i.e. a constant function when restricted to every Bernstein component of $\Irr(G)$). Then $\calZ_{comp}(G)=\calZ_{lc}(G)$.
\eth
\reft{main} has the following surprising corollary (in fact, technically we are first going to prove the corollary and then deduce \reft{main} from it, but historically our starting conjectural point was the assertion of \reft{main}):
\cor{mult}
$\calZ_{comp}(G)$ is a subalgebra of $\calZ(G)$.
\ecor

\refc{mult} is surprising since the set of compact elements of $G$ is not closed under multiplication.
We believe that \refc{mult} is actually a part of a more general statement.
While we are not sure what this statement really is, at least we believe in the following:
\conj{tnilp}
Let $\calZ_{tunip}(G)$ denote the subspace of $\calZ(G)$ consisting of distributions supported on topologically unipotent elements. Then
$\calZ_{tunip}(G)$ is a subalgebra of $\calZ(G)$.
\econj

\ssec{}{Relation to the work of J.-F.~Dat}
\reft{main} is in fact not completely new -- the inclusion $\calZ_{lc}(G)\subset \calZ_{comp}(G)$ was essentially proved by J.-F.~Dat
(cf. Section 2 of \cite{dat}).
 Namely, it is shown in {\em loc. cit.} that every idempotent in $\calZ(G)$ is supported on compact elements.
 Hence if for every $\Ome\in \calC(G)$ we denote by $\calE_{\Ome}$ the element of $\calZ(G)$ for which the function
 $f(\calE_{\Ome})$ is equal to 1 on $\Ome$ and is equal to 0 on any other component, then $\calE_{\Ome}\in \calZ_{comp}(G)$.
 On the other hand, any $\Phi\in \calZ(G)$ such that $f(\Phi)$ is constant on every $\Ome$ is locally on $G$ a linear combination of the distributions
 $\calE_{\Ome}$, hence $\Phi$ is supported on compact elements. The main observation of this note is that a mild adaptation of Dat's argument also proves the converse statement.

\ssec{}{A variant}In fact the inclusion $\calZ_{lc}(G)\subset \calZ_{comp}(G)$ has the following stronger version.
Given an element $g\in G$ we can define in a standard way a parabolic subgroup $P_g$ of $G$ and a strictly dominant element $\lambda_g\in Z(M_g)/Z(M_g)^0$ (the latter group is always a lattice and we shall denote the multiplication there by +; also, "strictly dominant" means
that the adjoint action of $\lam_g$ contract the unipotent radical of $P_g$ to the unit element).
Here we denote by $M_g$ the Levi group of $P_g$; also $Z(M)$ stands for the center of $M_g$ and $Z(M_g)^0$ is its maximal compact subgroup.
 Namely, $P_g$ consists of all $x\in G$ such that
$\lim\limits_{n\to \infty} g^n x g^{-n}$ exists. Also the image of $g$ under the natural map $P_g\to M_g$ must be compact modulo center  and
hence $g$ defines an element in $Z(M)/Z(M)\cap M^0=Z(M)/Z(M)^0$ which we call $\lam_g$. Note that $P_g=G$ if and only if $g$ is compact modulo center. Moreover, we have $P_g=G, \lam_g=0$ if and only if $g$ is compact.

Let now $\calP(G)$ denote the set of conjugacy classes of pairs $(P,\lam)$ as above. Then the above construction produces
a decomposition
\eq{decg}
G=\bigsqcup\limits_{(P,\lam)\in \calP(G)} G_{P,\lam}
\end{equation}
and each $G_{P,\lam}$ is an open subset of $G$ invariant under conjugation.

Let now $\calD(G)$ denote the space of distributions on $G$; let also $\calD^{inv}(G)\subset \calD(G)$ be the space of invariant distributions.
Then \refe{decg} produces a decomposition
\eq{decd}
\calD^{inv}(G)=\prod\limits_{(P,\lam)\in \calP(G)} \calD^{inv}_{P,\lam}.
\end{equation}
Here $\calD^{inv}_{P,\lam}$ consists of all invariant distributions supported on $G_{P,\lam}$.
We can now formulate
\th{main'}
Let $\Phi\in \calZ_{lc}(G)$. Then convolution with $\Phi$ preserves the decomposition \refe{decd} (i.e. preserves each $\calD^{inv}_{P,\lam}$).
\eth
This result is due to R.~Bezrukavnikov and we reproduce its proof in the Appendix.
\reft{main'} implies the inclusion $\calZ_{lc}(G)\subset\calZ_{comp}(G)$. Namely let $\del$ denote the delta-distribution at the unit element of $G$.
Obviously $\del\in \calD^{inv}_{G,0}$, hence $\Phi=\Phi\star\del\in \calD^{inv}_{G,0}$, i.e. $\Phi$ is supported on compact elements.
Our proof of \reft{main'} is somewhat simpler than the proof of the inclusion $\calZ_{lc}(G)\subset \calZ_{comp}(G)$ from \cite{dat}.
However, we still need the arguments of \cite{dat} in order to prove the opposite inclusion.
%------------------------------------------------------------------------------------------------------------------------
\ssec{}{An example}
For a rational number $r\in \QQ_{\geq 0}$ Moy and Prasad (cf. \cite{MP}) define a subset $\Irr_{\leq r}(G)$ of $\Irr(G)$ called "representations of depth $\leq r$". The set $\Irr_{\leq r}(G)$ is a union of components of $\Irr(G)$. Let $\Phi_r\in \calZ_G$ be the projector to $\Irr_{\leq r}(G)$; in other words $\Phi_r$ is the element of $\calZ(G)$ such that $f(\Phi_r)(\pi)=1$ if $\pi\in \Irr_{\leq r}(G)$ and $f(\Phi_r)(\pi)=0$ otherwise.
According to \reft{main} $\Phi_r$ should be concentrated on compact elements.
In \cite{BKV} the authors give an explicit formula for $\Phi_r$ which indeed shows this explicitly. In fact, the main result of \cite{BKV} implies a much stronger restriction to on the support of $\Phi_r$. It would be interesting to include this restriction into a general theorem in the style of \reft{main}.

\ssec{}{Geometric and spectral support}
We conclude the introduction with yet another conjecture which contains \reft{main} as a special case. To simplify the discussion we shall assume
that $G$ is a split.

Let us assume that $G=\bfG(\calK)$ where $\bfG$ is the corresponding split algebraic group defined over $\ZZ$.
Let $\Lam$ denote the coweight lattice of $\bfG$; we shall denote by $\Lam^+$ the set of dominant coweights. Also, for $\lam,\mu\in \Lam$
we shall write $\lam\geq\mu$ if $\lam-\mu$ is a sum of positive coroots of $\bfG$.
Let $K=\bfG(\calO)$. Then by Cartan decomposition the double quotient $K\backslash G/K$ is in natural bijection with the set
$\Lam^+$ of dominant coweights of $\bfG$. For each $\lam\in \Lam^+$ we shall denote by $G^{\lam}$ the corresponding double coset.
We set
$$
G^{\leq \lam}=\bigcup\limits_{\mu\in \Lam^+,\mu\leq \lam} G^{\mu}.
$$
We now define
\eq{geom}
\calZ_{geom}^{\leq \lambda}(G)=\{ \Phi\in\calZ(G)|\ \supp(\Phi)\subset Ad G\cdot (G^{\leq \lam})\}.
\end{equation}
Note that $\calZ_{geom}^{\leq 0}(G)=\calZ_{comp}(G)$.

On the hand, let $\bfP,\bfM$ be a parabolic subgroup of $G$ and its Levi subgroup (both defined over $\calK$).
Let $\Lam_M$ be the cocharacter lattice of $\bfM/[\bfM,\bfM]$. Then $\Lam_M$ is a sublattice of $\Lam$ and $\CC[\Lam_M]$ is the algebra of regular functions on the set of unramified characters of $M$ (for an element $\lam\in \Lam_M$ we shall denote by $e^{\lam}$ the corresponding element of $\CC[\Lam_M]$). Fix a cuspidal representation $\sig$ of $M$ with unitary central character. Then every $\Phi\in \calZ(G)$ defines an element $f_{\sig}(\Phi)\in \CC[\Lam_M]$. Namely, for an unramified $\psi:M\to \CC^*$ we define
$f_{\sig}(\Phi)(\psi)$ to be the scalar by which $\Phi$ acts in $i_{GM}(\sig\ten \psi)$ where $i_{GM}$ stands for unitary induction from $P$ to $G$.
We now define
\eq{spectral}
\calZ_{spectral}^{\leq \lam}(G)=\{ \Phi\in \calZ(G)| f_{\sig}(\Phi)=\sum\limits_{\mu\leq \lam} a_{\mu} e^{\mu}.
\end{equation}
It is easy to see that $\Phi\in \calZ_{spectral}^{\leq 0}(G)$ if and only if $f(\Phi)$ is constant on every component of $\Irr(G)$.
We can now formulate
\prop{spec-geom}
For any $\lam\in\Lam^+$ we have
$\calZ_{geom}^{\leq \lam}(G)=\calZ_{spectral}^{\leq \lam}(G)$.
\eprop
\refp{spec-geom} reduces to \reft{main} when $\lam=0$. The proof of \refp{spec-geom} easily follows from \reft{app}.
%------------------------------------------------------------------------------
\ssec{}{Acknowledgements}A.B. and R.B. were partially supported by the National Science Foundation. The project has received funding from ERC under grant agreement 669655. 
%---------------------------------------------------------------------------------------------------------------------
\sec{}{Proof of \refc{mult}}
 As was mentioned above, for the most part our proof is a repetition of the arguments of Section 2 of \cite{dat}.

\ssec{}{$\calZ_{comp}$ and $\ocalH$}
Let $\calH(G)$ denote the Hecke algebra of $G$. In what follows we shall choose a Haar measure on $G$ and we are going to identify $\calH(G)$ with the space of locally constant compactly supported functions on $G$.  Let $\ocalH(G)=\calH(G)/[\calH(G),\calH(G)]$. Obviously $\calZ(G)$ acts on $\ocalH$. For any $\Phi\in \calZ(G)$ we shall denote by $\oPhi$ the corresponding endomorphism of $\ocalH$.

Following \cite{dat} let us denote by $G^0$ the subgroup of $G$ generated by all the open compact subgroups of $G$. We also denote by
$G_c$ the set of compact elements modulo center. Then $G^0_c=G^0\cap G_c$ is the set of compact elements of $G$.

For an open subset $X$ of $G$ we denote by $1_X$ the characteristic function of $X$. Then multiplication (not convolution!) by $1_X$ is an endomorphism of $\calH$ which (abusing the notation) we shall also denote by the same symbol. Moreover, if $X$ is invariant under conjugation,
then multiplication by $1_X$ descends to endomorphism of $\ocalH$, which we shall denote by $\oone_X$. Then $\oone_{G^0}, \oone_{G_c}$ and $\oone_{G_c^0}$
are well-defined and we have $\oone_{G_c^0}=\oone_{G_c}\circ \oone_{G^0}$. Also, all of these endomorphisms commute with each other.

Below is the main result of this Section.

\th{main-comm}
We have $\Phi\in \calZ_{comp}(G)$ if and only if $\oPhi$ commutes with $\oone_{G_c^0}$.
\eth
It is clear that \reft{main-comm} implies \refc{mult}.

\ssec{}{Proof of \reft{main-comm}: the "if" direction}
This is the easy part of \reft{main-comm}. Note that two elements $h_1,h_2\in \calH$ have the same image in $\ocalH$ iff for any invariant distribution $\calE$ on $G$ we have $\calE(h_1)=\calE(h_2)$. In particular, if $\oh_1=\oh_2$, then $h_1(e)=h_2(e)$ where $e$ is the unit element of $G$.

Let us now assume that we are given $\Phi\in \calZ(G)$ such that $\oPhi$ commutes with $\oone_{G_c^0}$. Then for any $h\in \calH$ we have
$\overline{\Phi\star (h|_{G_c^0})}=\overline{(\Phi\star h)|_{G_c^0}}$. Let $\tilh(g)=h(g^{-1})$. Then we have
$$
\Phi(\tilh)=(\Phi\star h)(e)=(\Phi\star h)|_{G_c^0}(e)=(\Phi\star (h|_{G_c^0}))(e).
$$
Hence $\Phi(\tilh)=0$ if $\supp(\tilh)=\supp(h)\subset G\backslash G_c^0$, which means that $\supp \Phi\subset G_c^0$.
%----------------------------------------------------------------------------------------------------------------------
\ssec{easy}{}
Let us now start proving the opposite direction. Namely, let $\Phi\in \calZ_{comp}(G)$. We want to show that $\oPhi$ commutes with
$\oone_{G_c^0}$. For this it is enough to prove that $\oPhi$ commutes with $\oone_{G^0}$ and $\oone_{G_c}$. Let us first prove that
$\oPhi$ commutes with $\oone_{G^0}$. For this it is enough to prove that $\Phi$ commutes with $1_{G^0}$. In other words, we need to prove that
for any $h\in \calH$ we have $\Phi\star (h|_{G^0})=(\Phi\star h)|_{G^0}$. But this is obvious since $G^0$ is a subgroup of $G$ and $\supp(\Phi)\subset G^0$.
\ssec{}{Induction and restriction}Let us choose a split Cartan subgroup $T$ of $G$ and a Borel subgroup $B$ of $G$ with unipotent radical $U$. We denote by $\oB$ the opposite Borel subgroup of $G$.
Then we have a notion of standard Levi subgroup $M$ of $G$. We shall use the notation $M<G$ to indicate that $M$ is a standard Levi subgroup of $G$. To any such $M$ there corresponds a pair of parabolic subgroups $P,\oP$, where
$P\cap \oP=M$ and $B\subset P, \oB\subset \oP$. Also for any such $M$ we have maps (cf. \cite{dat} and references therein):
$i_{GM}^{\calZ},\oi_{GM}^{\calZ}:\calZ(M)\to \calZ(G),\  r_{GM}^{\calZ},\ovr_{GM}^{\calZ}:\calZ(G)\to \calZ(M), \ i_{GM}^{\ocalH},\oi_{GM}^{\ocalH}:\ocalH(M)\to\ocalH(G),\  r_{GM}^{\ocalH},\ovr_{GM}^{\ocalH}:\ocalH(G)\to\ocalH(M)$.
These maps satisfy the following properties:

\smallskip
0) All these maps are equal to identity when $M=G$.

1) For any $\Phi\in\calZ(G), f\in \ocalH(M)$ we have $\Phi\star \oi_{GM}^{\ocalH}(f)=\oi_{GM}^{\ocalH}(r_{GM}^{\calZ}(\Phi)\star f)$

2) For any $\Phi\in\calZ(G), h\in \ocalH(G)$ we have $r_{MG}^{\ocalH}(\Phi\star h)=r_{MG}^{\calZ}(\Phi)\star r_{MG}^{\ocalH}(h)$

3) Let $U_P$ denote the unipotent radical of a standard parabolic subgroup $P$ with the standard Levi subgroup $M$. Let $\pi_P$ denote the natural projection from $G$ to $G/U_P$; clearly $M$ is a closed subset of $G/U_P$. Let now $\Phi\in\calZ(G)$. The direct image $(\pi_P)_*\Phi$ makes sense - by the definition for any locally constant compactly supported function $\phi$ on $G/U_P$ we set $(\pi_P)_*\Phi(\phi)=(e_K\star \Phi)(\pi_P^*\phi)$,
where $K$ is any open compact subgroup of $G$ such that $\phi$ is $K$-invariant and $e_K$ is the Haar measure on $K$. Then the distribution
$(\pi_P)_*(\Phi)$ is concentrated on $M$. Moreover, the resulting distribution on $M$ is equal to $r_{GM}^{\calZ}(\Phi)$ up to multiplication by an unramified character of $M$.

\smallskip
\noindent
Property 3) above implies that if $\Phi\in\calZ_{comp}(G)$ then $r_{GM}^{\calZ}(\Phi)\in \calZ_{comp}(M)$. Indeed, this follows from the fact that
an element $g\in P\subset G$ is compact if and only if its projection to $M=P/U_P$ is compact.
%-----------------------------------------------------------------------------------------------------------
\ssec{}{Clozel's formula}
The main ingredient of the argument of Section 2 of \cite{dat} (and also of our proof of \reft{main-comm}) is the following formula due to Clozel.

\prop{clozel}
For any standard Levi $M$ there exists a function $\chi_M:M/M^0\to \CC$ such that for any $h\in \ocalH(G)$ we have
\eq{cloz}
h=\sum\limits_{M<G} \oi_{GM}^{\ocalH}(\chi_M\cdot (\oone_{M_c} (r_{GM}^{\ocalH}(h))).
\end{equation}
Moreover, $\chi_G=1$.
\eprop
\ssec{}{End of the proof}
We can now finish the proof of \reft{main-comm}. By induction we can assume that for any Levi subgroup $M$ of $G$ which is different from $G$ and any $\Phi\in \calZ_{comp}(M)$ we have $\oone_{M_c}(\Psi\star f)=\Psi\star (\oone_{M_c}(f))$ for any $f\in \ocalH(M)$. Let now $\Phi\in\calZ_{comp}(G), h\in \ocalH(G)$.
Then we have
$$
\begin{aligned}
\oone_{G_c}(\Phi\star h)=\Phi\star h - \sum\limits_{M<G, M\neq G} \oi_{GM}^{\ocalH}(\chi_M\cdot (\oone_{M_c} (r_{GM}^{\ocalH}(\Phi\star h)))=\\
\Phi\star h - \sum\limits_{M<G, M\neq G} \oi_{GM}^{\ocalH}(\chi_M\cdot (\oone_{M_c} (r_{GM}^{\calZ}(\Phi)\star r_{GM}^{\ocalH}(h))))\overset{1}=\\
\Phi\star h - \sum\limits_{M<G, M\neq G} \oi_{GM}^{\ocalH}(\chi_M\cdot (r_{GM}^{\calZ}(\Phi)\star\oone_{M_c} (r_{GM}^{\ocalH}(h)))\overset{2}=\\
\Phi\star h - \sum\limits_{M<G, M\neq G} \oi_{GM}^{\ocalH}(r_{GM}^{\calZ}(\Phi)\star \chi_M\cdot (\oone_{M_c} (r_{GM}^{\ocalH}(h)))=\\
\Phi\star h - \sum\limits_{M<G, M\neq G} \Phi\star \oi_{GM}^{\ocalH}(\chi_M\cdot (\oone_{M_c} (r_{GM}^{\ocalH}(h)))=\Phi\star (\oone_{G_c}(h)).
\end{aligned}
$$
Here $\overset{1}=$ follows from the fact that $r_{GM}^{\calZ}(\Phi)\in \calZ_{comp}(M)$ and from \refss{easy}. The equality $\overset{2}=$ follows from the induction hypothesis.
\sec{}{Proof of \reft{main}}
\ssec{}{}To finish the proof of \reft{main} we need to show that for any $\Phi\in \calZ_{comp}(G)$ the function $f(\Phi)$ is constant on every $\Ome\in\calC(G)$. Let us recall that we denote by $\calE_{\Ome}$ the element of $\calZ(G)$ for which the function
 $f(\calE_{\Ome})$ is equal to 1 on $\Ome$ and is equal to 0 on any other component. Then by \cite{dat} we have $\calE_{\Ome}\in\calZ_{comp}(G)$ and by \refc{mult} we also have $\Phi\star \calE_{\Ome}\in \calZ_{comp}(G)$ and it is enough to prove that $f(\Phi\star \calE_{\Ome})$ is constant on $\Ome$.
 Let $\calZ_{\Ome}(G)$ denote the subalgebra of $\calZ(G)$ consisting of elements $\Phi$ such that $f(\Phi)$ is equal to 0 on any $\Ome'\neq \Ome$ and let $\calZ_{\Ome,comp}(G)=\calZ_{comp}(G)\cap \calZ_{\Ome}(G)$. We already know that $\calE_{\Ome}\in \calZ_{\Ome,comp}(G)$. Clearly, our assertion follows from the following
\prop{omcomp}
$\calZ_{\Ome,comp}(G)=\CC\cdot \calE_{\Ome}$.
\eprop

We claim that \refp{omcomp} follows from the following:
\lem{fdim}
$\dim \calZ_{\Ome,comp}(G)< \infty$.
\elem
Indeed, $\calZ_{\Ome,comp}(G)$ is a subalgebra of $\calZ_{\Ome}(G)$ which is isomorphic to the algebra of functions on an irreducible algebraic variety $\bar{\Ome}$. Hence the only finite-dimensional subalgebra of it consists of constants. So, it remains to prove \refl{fdim}.
\prf
Let $\calD(G_c^0)$ denote the space of distributions on $G_c^0$. Consider the subspace $D_{\Ome}$ of $\calD(G_c^0)$ generated by distributions of the form
$\ch_{\pi}|_{G_c^0}$ where $\pi$ is an irreducible representation in $\Ome$ and $\ch_{\pi}$ is its character. Then it follows immediately from the Corollary in Section 3.1 in \cite{BDK} that $\calD_{\Ome}$ is finite-dimensional.

On the other hand, by Plancherel formula for every $\Ome\in \calC(G)$ there exists a measure $d\pi_{\Ome}$ on $\Ome$ such that
\eq{plancherel}
\del_G=\sum\limits_{\Ome\in \calC(G)}\int\limits_{\pi\in \Ome} \ch_{\pi} d\pi_{\Ome}.
\end{equation}
Here $\del_G$ denotes the $\delta$-distribution at the unit element of $G$.

Convolving this with a central element $\Phi\in \calZ(G)$ we get
\eq{plancherel-phi}
\Phi=\sum\limits_{\Ome\in \calC(G)}\int\limits_{\pi\in \Ome} f(\Phi)(\pi)\ch_{\pi} d\pi_{\Ome}.
\end{equation}
If $\Phi\in \calZ_{\Ome,comp}$ we get
\eq{plancherel-phi'}
\Phi=\int\limits_{\pi\in \Ome} f(\Phi)(\pi)\ch_{\pi} d\pi_{\Ome}.
\end{equation}
Since the LHS of \refe{plancherel-phi'} is concentrated on $G_c^0$, the same is true for RHS. Thus
we get
\eq{plancherel-phi''}
\Phi=\int\limits_{\pi\in \Ome} f(\Phi)(\pi)\ch_{\pi}|_{G_c^0} d\pi_{\Ome}.
\end{equation}
Hence $\Phi\in \calD_{\Ome}$, i.e. $\calZ_{\Ome,comp}\subset \calD_{\Ome}$, which implies that $\dim\calZ_{\Ome,comp}$ is finite-dimensional.
\epr
%---------------------------------------------------------------------------------------------------------
\sec{app}{Appendix: Proof of \reft{main'} (by R.~Bezrukavnikov)}
\ssec{}{Decomposition of $\ocalH(G)$}
We have an obvious perfect pairing between $\calD^{inv}(G)$ and $\ocalH(G)$. We claim that there is a decomposition
\eq{ohdec}
\ocalH(G)=\bigoplus\limits_{(P,\lam)\in \calP(G)} \ocalH(G)_{P,\lam},
\end{equation}
which is compatible with \refe{decd} by means of the above pairing. Namely, we let $\ocalH(G)_{P,\lam}$ to be the image of $\calH(G)_{P,\lam}$
where the latter consists of functions supported on $G_{P,\lam}$. The fact that \refe{ohdec} holds is clear.

\ssec{}{Spectral description of $\ocalH(G)$}
The space $\ocalH(G)$ admits the following well-known description. Let $\pi\in \calM(G)$ be a finitely generated representation and let $E$ be an endomorphism of $\pi$.
It is well-known (cf. e.g. \cite{Vig}) that we can associate to the pair $(\pi,E)$ and element $[\pi,E]$ of $\calH(G)$. Moreover, $\calH(G)$ is isomorphic to the $\CC$-span of symbols $[\pi,E]$ subject to the relations:

\smallskip
a) Let $\pi_1,\pi_2\in \calM(G)$ and let $u\in \Hom(\pi_1,\pi_2), v\in \Hom(\pi_2,\pi_1)$. Then $[\pi_1,vu]=[\pi_2,uv]$.

b) $[\pi_1,E_1]+[\pi_3,E_3]=[\pi_2,E_2]$ for a short exact sequence $0\to \pi_1\to \pi_2\to \pi_3\to 0$ which is compatible with the endomorphisms $E_i\in \End(\pi_i)$.

c) $[\pi, c_1E_1+c_2E_2]=c_1[\pi,E_1]+c_2[\pi,E_2]$, where $c_i\in \CC$ and $E_i\in \End(\pi)$.

\smallskip
\noindent
The action of $\calZ(G)$ on $\ocalH(G)$ can also be described in these terms. Namely, let $\Phi\in \calZ(G)$. Then
$\Phi\cdot [\pi,E]=[\pi,E\circ \pi(\Phi)]$.

In addition, let $\rho$ be an admissible representation of $G$. Then we have
\eq{char-cocenter}
\la [\pi,E],\ch_\rho\ra=\sum\limits_i (-1)^i \Tr(E,\Ext^i(\pi,\rho)).
\end{equation}

In view of the Trace Paley-Wiener theorem (cf. \cite{BDK}), \refe{char-cocenter} defines $[\pi,E]$ uniquely.
%----------------------------------------------------------------------------------------------------------------------------------------
\ssec{}{Spectral description of $\ocalH_{P,\lam}$}
Let now $P$, $\oP$ be a pair of opposite parabolic subgroups with $M=P\cap \oP$. Let $\lam\in Z(M)/Z(M)^0$ such
that $(\oP,\lam)\in \calP(G)$.
Let also $\sig$ be a finitely generated representation of $M$. Set
\eq{pi}
\pi=i_{GP}(\sig).
\end{equation}
Let us now choose a uniformizer $t$ of our local field. Then any $\lam\in Z(M)/Z(M^0)$ lifts naturally to an element $t^{\lam}\in Z(M)$. Hence it defines an endomorphism of $\sig$ and thus also of $\pi$. We shall denote this endomorphism by $E_{\lam}$.
\th{app}
The subspace $\ocalH_{\oP,\lam}$ is spanned by elements $[\pi,E_{\lam}]$ as above (here again $\oP$ denotes a parabolic subgroup which is opposite to $P$).
\eth

\noindent
{\bf Remark.} The element $[\pi,E_{\lam}]$ actually depends on the choice of $t$; however, it is easy to see that the span of all the $[\pi,E_{\lam}]$ does not.

\prf
For $(P,\lam)=(G,0)$ this is the "abstract Selberg principle" (cf. \cite{BlBr}).
The case $P=G$ and arbitrary $\lambda$ is completely analogous.

Let us now take arbitrary $P$ and $\lam$. Let $\sig$ be a finitely generated representation of the Levi group $M$ as above and $\lambda$ -- a strictly dominant cocharacter of $\pi$. Then we have a natural identification
$t^{\lam}M^0_c/ Ad(M)=G_{P,\lam}/Ad (G)$. Hence we get a natural isomorphism between $\ocalH(M)_{M,0}$ and $\ocalH(G)_{P,\lam}$. Indeed, if
an element of $\ocalH(M)$ is represented by some $h\in \calH(M)$ supported on $M^0_c=M_{M,0}$, then let us denote by $h_{\lam}$ the corresponding
element of $\calH(M)$ supported on $M_{M,\lam}=t^{\lam} \cdot M_c^0$.
For an open compact subgroup $K$ of $G$
let us denote by $h_{\lam,K}$ the result of averaging of $h_{\lam}$ with respect to the adjoint action of $K$. Its image in $\ocalH(G)$ is independent of $K$ and the assignment $h\text{ mod} [\calH(M),\calH(M)]\mapsto h_{\lam,K}\text{ mod} [\calH(G),\calH(G)]$ is the desired isomorphism.
Let us denote it by $\eta_{P,\lam}$.

Now in order to finish the proof it is enough to show that for $\pi$ as in \refe{pi} we have
\eq{brrr}
[\pi,E_{\lam}]=\eta_{\oP,\lam}([\sig,\Id]).
\end{equation}
Let $h=[\sig,\Id]$. Then it is easy to see that
\eq{brrrr}
[\sig,t^{\lam}]=h_{\lam}.
\end{equation}
Now to prove \refe{brrr} it is enough ( by Trace Paley-Wiener theorem) to check that both the LHS and the RHS of \refe{brrr} have the same inner product with
$\ch_{\rho}$ where $\rho$ stands for a generic irreducible representation of $G$. But we have
$$
\Ext^i_G(i_{GP}(\sig),\rho)=\Ext^i_M(\sig, r_{G\oP}(\rho)).
$$
Hence $\la [\pi,E_{\lam}],\ch_{\rho}\ra=\la [\sig,\lam],r_{G\oP}(\rho)\ra$ and \refe{brrr} follows from \refe{brrrr} and from the Casselman formula
for the character of  $r_{G\oP}(\rho)$ (cf. \cite{Cas}) which says for any $g\in G$ such that $P_g=\oP$ we have $\ch_{\rho}(g)=\ch_{r_{G\oP}(\rho)}(g)$.

\epr

\cor{}
\reft{main'} holds.
\ecor

\prf
It is enough to show that the action of any $\Phi\in\calZ_{lc}(G)$ preserves each $\ocalH_{\oP,\lam}$. Let us consider an element $[\pi,E]$ as above; without loss of generality we may assume that all irreducible subquotients of $\sig$ lie in one component of $\calC(M)$. But then all irreducible subquotients of $\pi$ as in lie in one component $\Ome\in \calC(G)$ and it follows that
$\Phi\star [\pi,E_{\lam}]=[\pi,f(\Phi)|_{\Ome}\cdot E_{\lam}]=f(\Phi)|_{\Ome}\cdot [\pi,E_{\lam}]$ (note that $f(\Phi)|_{\Ome}\in \CC$ as $\Phi\in \calZ_{lc}(G)$). Hence the span of all the $[\pi,E_{\lam}]$ is preserved by the convolution with $\Phi$.
\epr

\address{R.B.: Department of Mathematics, Massachusetts Institute of Technology}

\address{A.B.: Department of Mathematics,  University of Toronto, Perimeter Institute for Theoretical Physics and Department of Mathematics, Brown University}

\address{D.K.: Department of Mathematics, Hebrew University of Jerusalem}


\begin{thebibliography}{ddddd}

\bibitem{BD}
J.~Bernstein and P.~Deligne,
{\em Le ”centre” de Bernstein},
In ”Representations des groups reductifs sur un corps local,
Traveaux en cours” (P.Deligne ed.), Hermann, Paris, 1-32 (1984).

\bibitem{BDK}
J.~Bernstein, P.~Deligne and D.~Kazhdan,
{\em Trace Paley-Wiener theorem for reductive p-adic groups},
 J. Analyse Math. {\bf 47} (1986), 180-192.

\bibitem{BKV}
R.~Bezrukavnikov, D.~Kazhdan and Y.~Varshavsky,
{\em On the depth $r$ Bernstein projector}, arXiv:1504.01353.

\bibitem{BlBr}
P.~Blanc and J.-L.~Brylinsky,
{\em Cyclic homology and the Selberg principle},
J. Func. Anal. {\bf 109}, (1992) 289-330.

\bibitem{Cas}
W.~Casselman,
{\em Characters and Jacquet modules},
Math. Ann. {\bf 230} (1977), no. 2, 101-105.

\bibitem{dat}
J.-F.~Dat,
{\em Quelques propri\'et\'es des idempotents centraux
des groupes p-adiques},
J. reine angew. Math. {\bf 554} (2003), 69-103

\bibitem{MP}
A.~Moy and G.~Prasad,
{\em Unrefined minimal K-types for p-adic groups}, Invent. Math. {\bf 116}
(1994) 393–408.

\bibitem{Vig}
M.-F.~Vign\'eras,
{\em On formal dimensions for reductive p-adic groups},
Israel Math. Conf.
Proc. {\bf 2}, 225-265, 1990.
\end{thebibliography}
\end{document}